\documentclass{article}

\title  {Monoidal Morita invariants for finite group algebras}
\author {Kenichi Shimizu\footnote{Research Fellow of the Japan Society for the Promotion of Science}}
\date   {}

\usepackage{amsmath}
\usepackage{amssymb}
\usepackage{amsthm}
\usepackage{amscd}
\usepackage{graphicx}

\theoremstyle{plain}
\newtheorem{lemma}               {Lemma}[section]
\newtheorem{theorem}      [lemma]{Theorem}
\newtheorem{corollary}    [lemma]{Corollary}
\newtheorem{proposition}  [lemma]{Proposition}
\theoremstyle{definition}
\newtheorem{definition}   [lemma]{Definition}
\theoremstyle{remark}
\newtheorem{remark}       [lemma]{Remark}

\newcommand{\id}        {{\rm id}}
\newcommand{\op}        {{\rm op}}
\newcommand{\cop}       {{\rm cop}}
\newcommand{\Trace}     {{\rm Tr}}
\newcommand{\trace}     {{\rm tr}}
\newcommand{\bicross}   {{\mathop{\Join}\nolimits}}
\newcommand{\cotensor}  {{\mathop{\square}\nolimits}}
\newcommand{\End}       {\mathop\mathrm{End}\nolimits}
\newcommand{\Aut}       {\mathop\mathrm{Aut}\nolimits}
\newcommand{\Hom}       {\mathop\mathrm{Hom}\nolimits}
\newcommand{\Ext}       {\mathop\mathrm{Ext}\nolimits}
\newcommand{\Mod}       {\mathbf{Mod}}
\newcommand{\Com}       {\mathbf{Com}}
\newcommand{\finMod}    {\mathbf{mod}}

\begin{document}

\maketitle

\begin{abstract}
  Two Hopf algebras are called monoidally Morita equivalent if module categories over them are equivalent as linear monoidal categories. We introduce monoidal Morita invariants for finite-dimensional Hopf algebras based on certain braid group representations arising from the Drinfeld double construction. As an application, we show, for any integer $n$, the number of elements of order $n$ is a monoidal Morita invariant for finite group algebras. We also describe relations between our construction and invariants of closed 3-manifolds due to Reshetikhin and Turaev.
\end{abstract}

\section{Introduction}

Let $H$ be a Hopf algebra over a field $k$, for example, a group algebra. As is well known, the category of $H$-modules, denoted by ${\bf Mod}(H)$, is a $k$-linear monoidal category. We say that two Hopf algebras $H$ and $L$ are {\em monoidally Morita equivalent} if ${\bf Mod}(H)$ and ${\bf Mod}(L)$ are equivalent as $k$-linear monoidal categories. In this paper, we introduce monoidal Morita invariants of finite-dimensional Hopf algebras and apply them to finite group algebras.

Following Etingof and Gelaki \cite{MR1810480}, we say that two finite groups $G$ and $G'$ are {\em $k$-isocategorical} if $kG$ and $kG'$ are monoidally Morita equivalent. In the same paper, they classify all groups $\mathbb{C}$-isocategorical to a given finite group in group-theoretical terms. Following \cite{MR1810480}, we say that a finite group $G$ is {\em categorically rigid over $k$} if any group $k$-isocategorical to $G$ is isomorphic to $G$. As a direct consequence of their classification, $G$ is categorically rigid over $\mathbb{C}$ if $G$ does not admit a normal abelian subgroup $A$ of order $2^{2m}$ \cite[Corollary~1.4]{MR1810480}. In general, it is difficult to know when two finite groups are isocategorical even if we use the classification result. In this paper we show the following criterion, applying our results to finite group algebras.

\begin{theorem}
  \label{thm:main-theorem}
  Let $k$ be a field. If two finite groups $G$ and $G'$ are $k$-isocategorical, then for each positive integer $n$, the number of elements of order $n$ in $G$ is equal to the number of elements of order $n$ in $G'$.
\end{theorem}

Monoidal categories arise not only from algebra but also from low-di\-men\-sion\-al topology such as the theory of knots and braids. Our construction is based on certain braid group representations arising from the Drinfeld double $\mathcal{D}(H)$ of a finite-dimensional Hopf algebra $H$. Considering $\mathcal{D}(H)$ itself as a left $\mathcal{D}(H)$-module via the left multiplication, we have a series of canonical representations
\begin{equation*}
  \rho_n^{\mathcal{D}(H)}: B_n \to {\rm Aut}_{\mathcal{D}(H)}\left(\mathcal{D}(H)^{\otimes n}\right)
  \quad (n = 2, 3, \cdots)
\end{equation*}
of braid groups $B_n$. We show a monoidal Morita invariant $\tau(b; H)$ is given by $\tau(b; H) = {\rm Tr}(\rho_n^{\mathcal{D}(H)}(b))$, $b \in B_n$. Theorem~\ref{thm:main-theorem} is actually an application of these invariants associated with certain braids; see Section~\ref{sec:application}.

When $H$ is a finite-dimensional semisimple Hopf algebra over an algebraically closed field of characteristic zero, we can relate our construction to the Reshetikhin-Turaev invariant \cite{MR1292673} of closed 3-manifolds. This relation gives rise to the following theorem. For groups $X$ and $Y$, denote by ${\rm Hom}(X, Y)$ the set of group homomorphisms from $X$ to $Y$.

\begin{theorem}
  \label{thm:main-theorem-2}
  Let $k$ be a field. If finite groups $G$ and $G'$ are $k$-isocategorical, then for any oriented connected closed 3-manifold $M$, we have
  \begin{equation*}
    \# \Hom(\pi_1(M), G) = \# \Hom(\pi_1(M), G')
  \end{equation*}
  where $\pi_1(M)$ is the fundamental group of $M$.
\end{theorem}

This paper is organized as follows. In Section \ref{sec:preliminaries}, we introduce the notion of monoidal Morita invariance between Hopf algebras. We review Schauenburg's results \cite{MR1408508} and prove some lemmas for latter sections. In Section~\ref{sec:invariants}, we define monoidal Morita invariants associated with braids and introduce some basic properties of them. In Section~\ref{sec:application}, we apply our invariants to finite group algebras and prove Theorem~\ref{thm:main-theorem}. In Section~\ref{sec:topology}, we discuss relations between our invariants and the construction of invariants of closed 3-manifolds due to Reshetikhin and Turaev. The most part of this section is a review of \cite[Chapter~II]{MR1292673} and \cite[Chapter~4]{MR1797619}. Theorem~\ref{thm:main-theorem-2} will be proved in this section. Section~\ref{sec:examples} is devoted to further examples and applications of our invariants.

In Appendix~\ref{sec:similarity}, we argue similarity of permutation matrices and prove that two permutation matrices of same size are similar if and only if they are conjugate as permutations (Theorem~\ref{thm:similarity-theorem}). This theorem is used in Section~\ref{sec:application} as a part of the proof of Theorem~\ref{thm:main-theorem}.

Throughout this paper, the base field is denoted by $k$. Unless otherwise noted, vector spaces, algebras, coalgebras, etc. are over $k$. For vector spaces $V$ and $W$, $V \otimes W$ means $V \otimes_k W$. Functors between $k$-linear categories are always assumed to be $k$-linear. We use \cite{MR1786197} as a main reference for general theory of Hopf algebras. The comultiplication and counit of a bialgebra $H$ are denoted by $\Delta: H \to H \otimes H$ and $\varepsilon: H \to k$, respectively. The antipode of a Hopf algebra is denoted by $S$. We use Sweedler's sigma notation
\begin{equation*}
  \Delta(x) = \sum x_{(1)} \otimes x_{(2)}
\end{equation*}
to denote the comultiplication of an element $x$ in a coalgebra.

For an algebra $A$, we denote by $A^\op$ the opposite algebra. Similarly, for a coalgebra $C$, we denote by $C^\cop$ a coalgebra with the same underlying space with opposite comultiplication $\Delta^\cop$ given by $\Delta^\cop(c) = \sum c_{(2)} \otimes c_{(1)}$ for all $c \in C$ (the opposite coalgebra). For a bialgebra $H$, bialgebras $H^\op$ and $H^\cop$ are defined in an obvious way.

\section{Preliminaries}

\label{sec:preliminaries}

\subsection{Bialgebras and monoidal categories}

A {\em monoidal category} (or {\em tensor category}) is a category $\mathcal{C}$ equipped with a bifunctor $\otimes: \mathcal{C} \times \mathcal{C} \to \mathcal{C}$ and an object $\mathbf{1} \in \mathcal{C}$ satisfying certain associativity and unit constraints. The bifunctor $\otimes$ is called the {\em tensor product} and the object $\mathbf{1}$ is called the unit object. We refer the reader to Chapter XIII of Kassel \cite{MR1321145} for formal definitions of monoidal categories and monoidal functors.

The category of vector spaces, denoted by $\mathbf{Vec}(k)$, is a typical example of a $k$-linear monoidal category. A monoidal category is called {\em strict} if its associativity and unit isomorphisms are all identities. In this paper, we deal with mainly $k$-linear monoidal categories whose associativity and unit isomorphisms are ``trivial'' like $\mathbf{Vec}(k)$. Hence, although such categories are not strict, we treat them as if they were strict. This is valid since every monoidal category is equivalent to a strict one (see, e.g., \cite[XI.5]{MR1321145}).

Let $B$ be a bialgebra. Given left $B$-modules $V$ and $W$, the tensor product $V \otimes W$ is a left $B$-module by
\begin{equation*}
  x \cdot (v \otimes w) = \sum x_{(1)} v \otimes x_{(2)} w
\end{equation*}
for all $x \in B$, $v \in V$ and $w \in W$. The category of left $B$-modules, denoted by $\Mod(B)$, is a $k$-linear monoidal category with this tensor product. The unit object of $\Mod(B)$ is the trivial $B$-module $\mathbf{1} = k$ given by $x \cdot 1 = \varepsilon(x)1$ for all $x \in B$. The following proposition describes relations between monoidal categories $\Mod(B)$, $\Mod(B^\cop)$ and $\Mod(B^\op)$. For a monoidal category $\mathcal{C}$, we denote by $\mathcal{C}^{\rm rev}$ the monoidal category with the underlying category $\mathcal{C}$ and the reverse tensor product $\otimes^{\rm rev}$ given by $X \otimes^{\rm rev} Y = Y \otimes X$ for all $X, Y \in \mathcal{C}$.

\begin{proposition}
  \label{prop:mod-B-rev}
  Let $B$ be a bialgebra.

  {\rm (a)} $\Mod(B^\cop)$ is monoidally equivalent to $\Mod(B)^{\rm rev}$.

  {\rm (b)} $\Mod(B^\op)$  is monoidally equivalent to $\Mod(B)^{\rm rev}$ if $B$ has a bijective antipode.
\end{proposition}
\begin{proof}
  (a) The identity functor together with the monoidal structures
  \begin{equation*}
    T_{V,W}: V \otimes^{\rm rev} W = W \otimes V \to V \otimes W, \quad T_{V,W}(w \otimes v) = v \otimes w
  \end{equation*}
  gives a monoidal equivalence between $\Mod(B)$ and $\Mod(B^\cop)$.

  (b) Under this assumption, the antipode gives an isomorphism $B^\op \cong B^\cop$ of Hopf algebras. This induces a monoidal equivalence between $\Mod(B^\cop)$ and $\Mod(B^\op)$.
\end{proof}

We use the sigma notation such as $\rho(v) = \sum v_{(0)} \otimes v_{(1)}$ for right comodule structures. Given right $B$-comodules $V$ and $W$, the tensor product $V \otimes W$ is also a right $B$-comodule with structure map $\rho_{V \otimes W}^{}: V \otimes W \to V \otimes W \otimes B$ given by
\begin{equation*}
  \rho_{V \otimes W}^{}(v \otimes w) = \sum v_{(0)} \otimes w_{(0)} \otimes v_{(1)} w_{(1)}.
\end{equation*}
The category of right $B$-comodules is a monoidal category with this tensor product. We denote this monoidal category by $\Com(B)$. We can prove the following proposition in a similar way as Proposition \ref{prop:mod-B-rev}.

\begin{proposition}
  \label{prop:com-B-rev}
  Let $B$ be a bialgebra.

  {\rm (a)} $\Com(B^\op)$  is monoidally equivalent to $\Com(B)^{\rm rev}$.

  {\rm (b)} $\Com(B^\cop)$ is monoidally equivalent to $\Com(B)^{\rm rev}$ if $B$ has a bijective antipode.
\end{proposition}

\subsection{Monoidal Morita theory}

We introduce the following definition.

\begin{definition}
  Two Hopf algebras $H$ and $L$ are {\em monoidally Morita equivalent} if $\Mod(H)$ and $\Mod(L)$ are equivalent as ($k$-linear) monoidal categories.
\end{definition}

Although we are interested in modules over Hopf algebras, for a while, we refer to Schauenburg's results \cite{MR1408508} that deal with comodules over Hopf algebras. Let $C$ be a coalgebra. The {\em cotensor product} $V \cotensor_C^{} W$ of a right $C$-comodule $V$ and a left $C$-comodule $W$ is defined to be the kernel of
\begin{equation*}
  \rho_V^{} \otimes \id_W^{} - \id_V^{} \otimes \lambda_W^{}: V \otimes W \longrightarrow V \otimes C \otimes W
\end{equation*}
where $\rho_V^{}$ and $\lambda_W^{}$ are the structure maps of $V$ and $W$, respectively. Let $D$ be another coalgebra. If, moreover, $W$ is a $(C, D)$-bicomodule, the cotensor product $V \cotensor_C^{} W$ is naturally a right $D$-comodule.

Let $H$ be a Hopf algebra. For an $H$-comodule $M$, we let
\begin{equation*}
  M^{{\rm co}H} = \{ m \in M \mid \rho_M^{}(m) = m \otimes 1 \}
\end{equation*}
denote the space of $H$-coinvariants. A {\em right $H$-Galois object} is a right $H$-comodule algebra $A \ne 0$ such that $A^{{\rm co} H} = k$ and the linear map $A \otimes A \to A \otimes H$ given by $x \otimes y \mapsto \sum x y_{(0)} \otimes y_{(1)}$ is bijective. A {\em left  $H$-Galois object} is defined by replacing ``right'' with ``left''. For another Hopf algebra $L$, an $(H, L)$-biGalois object \cite[Definition~3.4]{MR1408508} is a left $H$- and right $L$-Galois object such that the two comodule structures make it an $(H, L)$-bicomodule. If $A$ is an $(H, L)$-biGalois object, the cotensor product functor
\begin{equation*}
  F_A: \Com(H) \to \Com(L), \quad F_A(V) = V \cotensor_H^{} A
\end{equation*}
gives a monoidal equivalence together with the monoidal structures
\begin{equation*}
  \begin{split}
    J_{V,W}^{}: (V \cotensor_H^{} A) \otimes (W \cotensor_H^{} A)
    & \to (V \otimes W) \cotensor_H^{} A, \\
    \left( \sum v_i \otimes x_i \right) \otimes \left( \sum w_j \otimes y_j \right)
    & \mapsto \sum v_i \otimes w_j \otimes x_i y_j
  \end{split}
\end{equation*}
and $\varphi^{}: k \to k \cotensor_H^{} A$, $a \mapsto a \otimes 1$.

\begin{theorem}[Schauenburg {\cite[Corollary~5.7]{MR1408508}}]
  \label{thm:Schauenburg}
  The above correspondence $A \mapsto (F_A, J, \varphi)$ gives a bijection between isomorphism classes of $(H, L)$-biGalois objects and isomorphism classes of monoidal equivalences $\Com(H) \to \Com(L)$.
\end{theorem}

A right $H$-Galois object $A$ is said to be {\em cleft} if there exists a convolution invertible $H$-colinear map $H \to A$ where we consider $H$ as a right $H$-comodule via the comultiplication. Note that $A$ is cleft if and only if $A$ is isomorphic to $H$ as a right $H$-comodule \cite[Theorem~9]{MR834465}. The notion of {\em cleft left $H$-Galois objects} is defined similarly. If $H$ is finite-dimensional, all left $H$-Galois objects and all right $H$-Galois objects are cleft. In the following lemma, we denote by $U_H$ the forgetful functor $\Com(H) \to \mathbf{Vec}(k)$.

\begin{lemma}
  \label{lem:key-lemma-com}
  Let $H$ and $L$ be finite-dimensional Hopf algebras. For any monoidal equivalence $F: \Com(H) \to \Com(L)$, the followings hold.

  {\rm (a)} $F(H)$ is isomorphic to $L$ in $\Com(L)$.

  {\rm (b)} $U_L \circ F$ is isomorphic to $U_H$ as a $k$-linear functor.
\end{lemma}
\begin{proof}
  By Theorem~\ref{thm:Schauenburg}, there exists an $(H, L)$-biGalois object $A$ such that $F$ is isomorphic to $F_A = (-) \cotensor_H^{} A$. Since $H$ and $L$ are finite-dimensional, $A$ is cleft. In particular, $A$ is isomorphic to $H$ as a left $H$-comodule and is isomorphic to $L$ as a right $L$-comodule.

  (a) We have $F(H) \cong H \cotensor_H^{} A \cong A \cong L$ as right $L$-comodules.

  (b) It can be proved in a similar way as in (a).
\end{proof}

We now return to modules over Hopf algebras. Let $H$ be a finite-dimensional Hopf algebra. Recall that we can identify $\Mod(H)$ with $\Com(H^*)$ where $H^*$ is the dual Hopf algebra of $H$. In particular, we consider a right $H^*$-comodule $V$ as a left $H$-module by
\begin{equation*}
  x \cdot v = \sum v_{(0)} \langle v_{(1)}, x \rangle
\end{equation*}
for all $x \in H$ and $v \in V$. In the following lemma, we denote by $U'_H: \Mod(H) \to \mathbf{Vec}(k)$ the forgetful functor.

\begin{lemma}
  \label{lem:key-lemma}
  Let $H$ and $L$ be finite-dimensional Hopf algebras. For any monoidal equivalence $F: \Mod(H) \to \Mod(L)$, the followings hold.

  {\rm (a)} $F(H)$ is isomorphic to $L$ in $\Mod(L)$.

  {\rm (b)} $U'_L \circ F$ is isomorphic to $U'_H$ as a $k$-linear functor.
\end{lemma}
\begin{proof}
  (a) Let $F: \Mod(H) \to \Mod(L)$ be an equivalence of monoidal categories. If we identify $\Mod(H)$ and $\Mod(L)$ with $\Com(H^*)$ and $\Com(L^*)$ respectively, we have an equivalence $F: \Com(H^*) \to \Com(L^*)$ of monoidal categories. By Lemma~\ref{lem:key-lemma-com} (a), $F(H^*) \cong L^*$ in $\Com(L^*)$. Since $H^* \in \Com(H^*)$ is isomorphic to $H$ as a left $H$-module (see \cite[Chapter 5]{MR1786197}), we have $F(H) \cong F(H^*) \cong L^* \cong L$ as left $L$-modules.

  (b) This is obvious by Lemma~\ref{lem:key-lemma-com} (b).
\end{proof}

\subsection{Braiding}

A {\em braiding} in a (strict) monoidal category $\mathcal{C}$ is a natural isomorphism $c_{X,Y}^{}: X \otimes Y \to Y \otimes X$ ($X, Y \in \mathcal{C}$) satisfying equations $c_{X \otimes Y, Z}^{} = (c_{X,Z}^{} \otimes \id_Y^{}) (\id_X^{} \otimes c_{Y,Z}^{})$ and $c_{X, Y \otimes Z}^{} = (\id_Y^{} \otimes c_{X,Z}^{}) (c_{X,Y}^{} \otimes \id_Z^{})$ for $X, Y, Z \in \mathcal{C}$. A {\em braided monoidal category} (or {\em braided tensor category}) is a monoidal category equipped with a braiding (see \cite[Chapter~XIII]{MR1321145}).

Let $\mathcal{C}$ be a braided monoidal category with braiding $c$. Each object of $\mathcal{C}$ yields a series of representations of braid groups. We denote by $B_n$ ($n \ge 2$) the braid group on $n$ strands. As is well known, $B_n$ is generated by basic braids
\begin{equation*}
  \sigma_i = \begin{array}{c} \mbox{\includegraphics{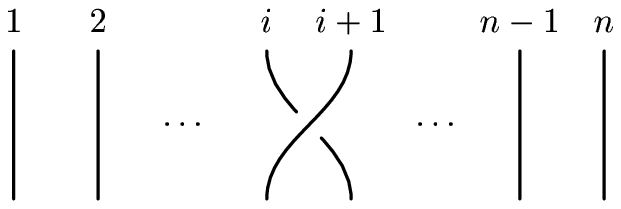}} \end{array}
  \quad (i = 1, 2, \cdots, n-1)
\end{equation*}
with defining relations
\begin{equation*}
  \begin{array}{cc}
    \sigma_i \sigma_j = \sigma_j \sigma_i & \quad \mbox{(if $|i-j|>1$)}, \\
    \sigma_i \sigma_{i+1} \sigma_i = \sigma_{i+1} \sigma_i \sigma_{i+1} & \quad \mbox{($i = 1, 2, \cdots, n-2$)}.
  \end{array}
\end{equation*}
Fix an object $X \in \mathcal{C}$. Then the morphism $\sigma = c_{X,X}^{}$ is a solution of the Yang-Baxter equation $(\sigma \otimes \id_X) (\id_X \otimes \sigma) (\sigma \otimes \id_X) = (\id_X \otimes \sigma) (\sigma \otimes \id_X) (\id_X \otimes \sigma)$ in $\Aut_\mathcal{C}(X^{\otimes 3})$, and thus we have a series of representations
\begin{equation*}
  \rho_n^X: B_n \to \Aut_{\mathcal{C}}(X^{\otimes n}) \quad (n = 2, 3, \cdots)
\end{equation*}
given by
\begin{equation*}
  \rho_n^{X} (\sigma_i) = \underbrace{\id_X \otimes \cdots \otimes \id_X}_{i-1}
  \otimes \, c_{X,X} \otimes \underbrace{\id_X \otimes \cdots \otimes \id_X}_{n-i-1}.
\end{equation*}

\begin{remark}
  \label{rem:braid-group-repr-1}
  Let $f: X \to Y$ be a morphism in $\mathcal{C}$. Since the braiding is natural, the diagram
  \begin{equation*}
    \begin{CD}
      X^{\otimes n} @>{\rho_n^X(b)}>> X^{\otimes n} \\
      @V{f \otimes \cdots \otimes f}VV @VV{f \otimes \cdots \otimes f}V \\
      Y^{\otimes n} @>>{\rho_n^Y(b)}> Y^{\otimes n} \\
    \end{CD}
  \end{equation*}
  commutes for all $b \in B_n$. In particular, braid group representations $\rho_n^X$ and $\rho_n^Y$ are equivalent if $X$ and $Y$ are isomorphic.
\end{remark}

\begin{remark}
  \label{rem:braid-group-repr-2}
  Let $\mathcal{D}$ be another braided monoidal category and $F: \mathcal{C} \to \mathcal{D}$ be a braided monoidal functor. If we fix an object $X \in \mathcal{C}$, we have representations
  \begin{equation*}
    \rho'_n: B_n \to \Aut_{\mathcal{D}}(F(X)^{\otimes n}),
    \quad \rho'_n(b) = \rho_n^{F(X)}(b)
  \end{equation*}
  and
  \begin{equation*}
    \rho''_n: B_n \to \Aut_{\mathcal{D}}(F(X^{\otimes n})),
    \quad \rho''_n(b) = F(\rho_n^{X}(b)).
  \end{equation*}
  These representations are equivalent. In fact, the canonical isomorphism $F(X)^{\otimes n} \cong F(X^{\otimes n})$ given by the monoidal structures of $F$ gives an intertwiner.
\end{remark}

\subsection{Quasitriangular Hopf algebras}

We argue braidings in the category of modules over a Hopf algebra. If $R = \sum s_i \otimes t_i \in A^{\otimes 2}$ is a universal $R$-matrix \cite[Definition~VIII.2.2]{MR1321145} of $A$, $\Mod(A)$ is a braided monoidal category with braiding $c^R_{V,W}: V \otimes W \to W \otimes V$ given by
\begin{equation*}
  c^R_{V,W}(v \otimes w) = \sum t_i w \otimes s_i v
\end{equation*}
for $v \in V$ and $w \in W$. It is known that this gives a one-to-one correspondence between braidings of $\Mod(A)$ and universal $R$-matrices of $A$. We denote this braided monoidal category by $\Mod(A, R)$. We often omit $R$ and denote $\Mod(A, R)$ by $\Mod(A)$ if $R$ is obvious.

A {\em quasitriangular Hopf algebra} is a pair $(A, R)$ of a Hopf algebra $A$ and a universal $R$-matrix of $A$. We list basic properties of quasitriangular Hopf algebras.

\begin{proposition}[{\cite[\S 2]{MR1025154}}, {\cite[Chapter VIII]{MR1321145}}]
  \label{prop:qt-hopf}
  Let $(A, R)$ be a quasitriangular Hopf algebra. Set $u = \sum S(t_i) s_i$ where $\sum s_i \otimes t_i = R$. Then the followings hold.

  {\rm (a)} $u$ is invertible with inverse $u^{-1} = \sum t_i S^2(s_i)$.

  {\rm (b)} The antipode $S$ is bijective and we have $S^2(x) = u x u^{-1}$ for all $x \in A$.

  {\rm (c)} $(\varepsilon \otimes \id_A)(R) = 1 = (\id_A \otimes \varepsilon)(R)$.

  {\rm (d)} $(S \otimes \id_A)(R) = R^{-1} = (\id_A \otimes S^{-1})(R)$ and $(S \otimes S)(R) = R$.
\end{proposition}

The element $u$ above is called the {\em Drinfeld element} of $(A, R)$. Note that $A$ is {\em involutive}, i.e., $S^2 = \id_A$ if and only if $u$ is central.

Let $\mathcal{C}$ be a braided monoidal category with braiding $c$. Then the reverse monoidal category $\mathcal{C}^{\rm rev}$ is also a braided monoidal category with braiding $c^{\rm rev}_{X,Y} = c_{Y,X}^{}: X \otimes^{\rm rev} Y \to Y \otimes^{\rm rev} X$. $\mathcal{C}^{\rm rev}$ is equivalent to $\mathcal{C}$ as a braided monoidal category. In fact, the identity functor together with monoidal structure $\varphi_{V,W}^{} = c_{W,V}^{}: V \otimes^{\rm rev} W \to V \otimes W$ gives an equivalence.

\begin{proposition}
  \label{prop:braided-equiv}
  Let $(A, R)$ be a quasitriangular Hopf algebra.

  {\rm (a)} $(A^{\cop}, R_{21})$ is isomorphic to $(A^{\op}, R_{21})$ as a quasitriangular Hopf algebra.

  {\rm (b)} $\Mod(A, R)$, $\Mod(A^{\cop}, R_{21})$ and $\Mod(A^{\op}, R_{21})$ are equivalent as braided monoidal categories.
\end{proposition}
\begin{proof}
  (a) The antipode $S: A^{\cop} \to A^{\op}$ gives an isomorphism of Hopf algebras. This preserves the universal $R$-matrix since $(S \otimes S)(R) = R$.

  (b) The equivalence given in the proof of Proposition~\ref{prop:mod-B-rev} induces an equivalence between braided monoidal categories $\Mod(A^{\cop}, R_{21})$ and $\Mod(A, R)^{\rm rev}$. The latter is equivalent to $\Mod(A, R)$ as we remarked above.
\end{proof}

\begin{lemma}
  \label{lem:braid-group-repr-3}
  Let $(A, R)$ and $(A', R')$ be finite-dimensional quasitriangular Hopf algebras. If $\Mod(A, R)$ and $\Mod(A', R')$ are equivalent as braided monoidal categories, braid group representations $\rho_n^A$ and $\rho_n^{A'}$ are equivalent for each $n \ge 2$.
\end{lemma}
\begin{proof}
  Let $F: \Mod(A, R) \to \Mod(A', R')$ be an equivalence of braided monoidal categories. By Lemma~\ref{lem:key-lemma}, Remark~\ref{rem:braid-group-repr-1} and Remark~\ref{rem:braid-group-repr-2}, we have isomorphisms $\eta_n: A^{\otimes n} \to F(A^{\otimes n})$, $\eta'_n: F(A^{\otimes n}) \to F(A)^{\otimes n}$ and $\eta''_n: F(A)^{\otimes n} \to A'^{\otimes n}$ such that the diagram in $\mathbf{Vec}(k)$
  \begin{equation*}
    \begin{CD}
      A^{\otimes n} @>{\eta_n}>> F(A^{\otimes n}) @>{\eta'_n}>> F(A)^{\otimes n} @>{\eta''_n}>> A'^{\otimes n} \\
      @V{\rho_n^A(b)}VV @V{F(\rho_n^A(b))}VV @VV{\rho_n^{F(A)}(b)}V @VV{\rho_n^{A'}(b)}V \\
      A^{\otimes n} @>>{\eta_n}> F(A^{\otimes n}) @>>{\eta'_n}> F(A)^{\otimes n} @>>{\eta''_n}> A'^{\otimes n}
    \end{CD}
  \end{equation*}
  commutes for all $b \in B_n$.
\end{proof}

\section{Invariants associated with braids}

\label{sec:invariants}

In this section, we define monoidal Morita invariants of finite-dimensional Hopf algebras associated with braids. Our construction is based on braid group representations arising from quasitriangular structures. A Hopf algebra does not always have universal $R$-matrices. We recall the Drinfeld double construction \cite[Chapter~IX]{MR1321145} which admits the canonical quasitriangular structure.

For a finite-dimensional Hopf algebra $H$, let $\mathcal{D}(H)$ be the Drinfeld double of $H$. Recall that $\mathcal{D}(H) = H^{*\cop} \otimes H$ as a coalgebra. To avoid confusion, we denote $f \otimes x \in \mathcal{D}(H)$ by $f \bicross x$. $\mathcal{D}(H)$ has a universal $R$-matrix
\begin{equation*}
  \mathcal{R}(H)
  = \sum_{i=1}^n \varepsilon \bicross h_i \otimes h_i^* \bicross 1 \in \mathcal{D}(H) \otimes \mathcal{D}(H)
\end{equation*}
where $\{ h_1, \cdots, h_n \}$ is a basis of $H$ and $\{ h_1^*, \cdots, h_n^* \}$ is the dual basis. $\mathcal{R}(H)$ is denoted by $\mathcal{R}$ if $H$ is obvious. Note that the braided monoidal category $\Mod(\mathcal{D}(H), \mathcal{R})$ is characterized as the {\em center}, denoted by $\mathcal{Z}(\Mod(H))$, of the monoidal category $\Mod(H)$ \cite[XIII.5]{MR1321145}. If finite-dimensional Hopf algebras $H$ and $L$ are monoidally Morita equivalent, we have equivalences
\begin{equation*}
  \Mod(\mathcal{D}(H), \mathcal{R}) \approx \mathcal{Z}(\Mod(H))
  \approx \mathcal{Z}(\Mod(L)) \approx \Mod(\mathcal{D}(L), \mathcal{R})
\end{equation*}
of braided monoidal categories. Applying Lemma~\ref{lem:braid-group-repr-3}, we have the following theorem.

\begin{theorem}
  Let $H$ and $L$ be finite-dimensional Hopf algebras. If $H$ and $L$ are monoidally Morita equivalent, then, for any integers $n \ge 2$, braid group representations
  \begin{equation*}
    \rho_n^{\mathcal{D}(H)}: B_n \to \Aut(\mathcal{D}(H)^{\otimes n}) \quad \mbox{and} \quad
    \rho_n^{\mathcal{D}(L)}: B_n \to \Aut(\mathcal{D}(L)^{\otimes n})
  \end{equation*}
  are equivalent.
\end{theorem}

\begin{remark}
  The above theorem gives us various monoidal Morita invariants. For instance, the exponent of finite-dimensional Hopf algebra $H$, which is defined to be the smallest integer $n \ge 1$ such that the equation
  \begin{equation*}
    \sum x_{(1)} S^{-2}(x_{(2)}) \cdots S^{-2 n + 2}(x_{(n)}) = \varepsilon(x) 1_H
  \end{equation*}
  holds for all $x \in H$, is equal to the order of $\rho_n^{\mathcal{D}(H)}(\sigma_1^2)$ (see \cite[Theorem 2.5]{MR1689203}).
\end{remark}
We study the following type of invariants.

\begin{definition}
  Let $b \in B_n$ be a braid. We define a {\em monoidal Morita invariant $\tau(b; H)$ associated with $b$ and a finite-dimensional Hopf algebra $H$} by
  \begin{equation*}
    \tau(b; H) = \Trace(\rho_n^{\mathcal{D}(H)}(b)).
  \end{equation*}
\end{definition}

Let us list some elementary properties of $\tau(b; H)$. For braids $b_1 \in B_n$ and $b_2 \in B_m$, we denote by $b_1 \otimes b_2 \in B_{n+m}$ the braid on $n+m$ strands which is obtained by arranging $b_2$ to the right of $b_1$.

\begin{proposition}
  \label{prop:basic-properties}
  Let $H$ be a finite-dimensional Hopf algebra. Then:

  {\rm (a)} $\tau(1_n; H) = \dim(H)^{2n}$ where $1_n$ is the identity of $B_n$.

  {\rm (b)} $\tau(b_1 b_2; H) = \tau(b_2 b_1; H)$ for all $b_1, b_2 \in B_n$.

  {\rm (c)} $\tau(b_1 \otimes b_2; H) = \tau(b_1; H) \tau(b_2; H)$ for all $b_1 \in B_n$ and $b_2 \in B_m$.

  {\rm (d)} If $K/k$ is a field extension, $\tau(b; K \otimes_k H) = \tau(b; H)$.
\end{proposition}
\begin{proof}
  Proofs are obvious from properties of trace.
\end{proof}

Our invariants cannot distinguish a finite-dimensional Hopf algebra and its dual since the construction is based on the Drinfeld double.

\begin{proposition}
  Let $H$ be a finite-dimensional Hopf algebra. $\Mod(\mathcal{D}(H), \mathcal{R})$ and $\Mod(\mathcal{D}(H^*), \mathcal{R})$ are equivalent as braided monoidal categories.
\end{proposition}
\begin{proof}
  Recall that $\mathcal{D}(H) = H^* \otimes H$ as a vector space. Under the canonical identification $H^{**} \cong H$, a linear map $T: H^* \otimes H \to H \otimes H^*$ given by $T(f \otimes x) = x \otimes f$ induces an isomorphism
  \begin{equation*}
    T: (\mathcal{D}(H), \mathcal{R}(H)) \to (\mathcal{D}(H^{\op \cop *})^{\op}, \mathcal{R}(H^{\op \cop *})_{21})
  \end{equation*}
  of quasitriangular Hopf algebras \cite[Theorem~3]{MR1220770}. Since $H$ is finite-dimensional, $H$ and $H^{\op \cop}$ are isomorphic as Hopf algebras via the antipode. Therefore, we have an isomorphism $(\mathcal{D}(H^*), \mathcal{R}(H^*)) \cong (\mathcal{D}(H)^{\op}, \mathcal{R}(H)_{21})$ of quasitriangular Hopf algebras. Applying Proposition~\ref{prop:braided-equiv} completes the proof.
\end{proof}

\begin{corollary}
  For each $n \ge 2$, braid group representations $\rho_n^{\mathcal{D}(H)}$ and $\rho_n^{\mathcal{D}(H^*)}$ are equivalent.
  In particular, $\tau(b; H^*) = \tau(b; H)$ for any braid $b$.
\end{corollary}

\section{Application to group algebras}

\label{sec:application}

\begin{figure}[htbp]
  \centering
  \includegraphics{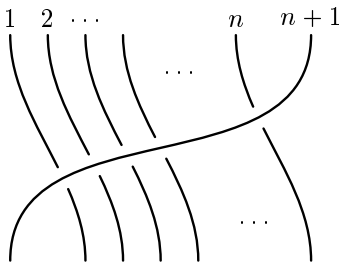}
  \caption{The braid $b_n \in B_{n+1}$}
  \label{fig:lens-braid}
\end{figure}

Our aim in this section is to prove Theorem \ref{thm:main-theorem}. This will be done by calculating the monoidal Morita invariant associated with braid
\begin{equation*}
  b_n = \sigma_n \sigma_{n-1} \cdots \sigma_1 \in B_{n+1} \quad (n = 1, 2, \cdots)
\end{equation*}
which is illustrated as Figure~\ref{fig:lens-braid}.

\subsection{Reduction to characteristic zero}

Let $G$ be a finite group. As we prove later in Lemma~\ref{lem:trace-3}, we have
\begin{equation*}
  \tau(b_n; kG) = |G| \cdot \# \{ g \in G \mid g^n = 1 \}
\end{equation*}
for each $n$. This equation holds in any characteristic. In characteristic zero, Theorem~\ref{thm:main-theorem} follows easily from it. The first step of the proof of Theorem~\ref{thm:main-theorem} is to reduce the problem to the case when the characteristic of $k$ is zero. We have the following theorem.

\begin{theorem}
  \label{thm:base-change-theorem}
  {\rm (a)} Let $G$ be a finite group. Then, $\tau(b; \mathbb{C}G)$ is a non-negative integer for any braid $b \in B_n$. For an arbitrary field $k$, we have $\tau(b; \mathbb{C}G) = \tau(b; kG)$ in $k$.

  {\rm (b)} Let $k$ be an arbitrary field. If two finite groups $G$ and $G'$ are $k$-isocategorical, we have $\tau(b; \mathbb{C}G) = \tau(b; \mathbb{C}G')$ for any braid $b \in B_n$.
\end{theorem}

In fact, the field $\mathbb{C}$ in Theorem~\ref{thm:base-change-theorem} (b) can be replaced by any field. There are two reasons why we use $\mathbb{C}$. First, $\mathbb{C}$ is an algebraically closed field of characteristic zero. Second, we desire to relate our monoidal Morita invariants to certain theories of closed 3-manifolds; see Section~\ref{sec:topology}.

We recall the structure of $\mathcal{D}(kG)$. For each $g \in G$, define $e_g \in (kG)^*$ by $\langle e_g, x \rangle = \delta_{x, g}$ for all $x \in G$ where $\delta$ is the Kronecker delta. Then, the set $\{ e_g \bicross x \}_{g, x \in G}$ is a basis of $\mathcal{D}(kG)$. The multiplication of $\mathcal{D}(kG)$ is given by
\[ (e_g \bicross x) (e_h \bicross y) = \delta_{g, xhx^{-1}} \, e_g \bicross (xy) \]
for all $g, h, x, y \in G$. The comultiplication $\Delta$ is given by
\[ \Delta(e_g \bicross x) = \sum_{h \in G} (e_{hg} \bicross x) \otimes (e_{gh^{-1}} \bicross x) \]
for all $g, x \in G$. Set $1^* = \sum_{g \in G} e_g$ (this is the counit of $kG$). Then, the universal $R$-matrix of $\mathcal{D}(kG)$ is given by
\[ \mathcal{R}(kG) = \sum_{g \in G} (1^* \bicross g) \otimes (e_g \bicross 1). \]
The Drinfeld element $u$ and its inverse are given respectively by
\[ u = \sum_{g \in G} e_g \bicross g^{-1} \quad \mbox{and} \quad u^{-1} = \sum_{g \in G} e_g \bicross g. \]

The proof of Theorem~\ref{thm:base-change-theorem} is based on the following observation: For all braid $b \in B_n$, $\rho_n^{\mathcal{D}(kG)}(b)$ is represented by a permutation matrix in basis
\[ e_{g_1} \bicross x_1 \otimes \cdots \otimes e_{g_n} \bicross x_n \quad (g_i, x_i \in G). \]
Note that this permutation is independent from the base field $k$. For a permutation matrix $P$, we denote by ${\rm Fix}(P)$ the number of fixed points of the corresponding permutation. The following lemma is a direct consequence of Theorem~\ref{thm:similarity-theorem} in Appendix A.

\begin{lemma}
  \label{lem:fixed-points}
  Let $P$ and $Q$ be permutation matrices of the same size. If $P$ and $Q$ are similar over $k$, we have ${\rm Fix}(P) = {\rm Fix}(Q)$.
\end{lemma}

If the characteristic of the base field $k$ is zero, the proof is obvious since $\Trace(P) = {\rm Fix}(P)$ in $k$. Lemma~\ref{lem:fixed-points} allows us to define ${\rm Fix}(P)$ for an automorphism $P$ on a finite-dimensional vector space which is represented by a permutation matrix in some basis. Now we can prove Theorem~\ref{thm:base-change-theorem}.

\begin{proof}[Proof of Theorem~\ref{thm:base-change-theorem}]
  (a) Let $b \in B_n$ be a braid. $\tau(b; \mathbb{C}G)$ is a non-negative integer as the trace of a permutation matrix. Since, as we remarked above, a permutation induced by $\rho_n^{\mathcal{D}(kG)}(b)$ is independent from the choice of the base field $k$, we have
  \begin{equation*}
    {\rm Fix}\left(\rho_n^{\mathcal{D}(kG)}(b)\right)
    = {\rm Fix}\left(\rho_n^{\mathcal{D}(\mathbb{C}G)}(b)\right)
    = \Trace\left(\rho_n^{\mathcal{D}(\mathbb{C}G)}(b)\right)
    = \tau(b; \mathbb{C}G).
  \end{equation*}
  Thus, we have $\tau(b; kG) = \tau(b; \mathbb{C}G)$ in $k$.

  (b) Assume that finite groups $G$ and $G'$ are $k$-isocategorical. Theorem~\ref{prop:basic-properties} yields
  \begin{equation*}
    {\rm Fix}\left(\rho_n^{\mathcal{D}(kG)}(b)\right) = {\rm Fix}\left(\rho_n^{\mathcal{D}(kG')}(b)\right)
  \end{equation*}
  for all $b \in B_n$. Thus, we have $\tau(b; \mathbb{C}G) = \tau(b; \mathbb{C}G')$.
\end{proof}

\subsection{The number of elements of order $n$}

Let $(A, R)$ be a quasitriangular Hopf algebra and $V$ be a finite-dimensional left $A$-module. Write $R = \sum_i s_i \otimes t_i$. By the definition of $\rho_{n+1}^V$, we have
\begin{equation*}
  \rho_{n+1}^V(b_n)(v_0 \otimes \cdots \otimes v_n)
  = \sum_{i_1, \cdots, i_n} t_{i_1} v_1 \otimes \cdots \otimes t_{i_n} v_n \otimes s_{i_n} \cdots s_{i_1} v_0
\end{equation*}
for all $v_0, \cdots, v_n \in V$. First, we give a description of the trace of $\rho_{n+1}^V(b_n)$ in terms of $R$ and the Drinfeld element $u$. For $a \in A$, we denote by $\Trace_V(a)$ the trace of the linear endomorphism on $V$ given by $v \mapsto a \cdot v$ ($v \in V$).

\begin{lemma}
  \label{lem:trace-1}
  Notations are as above.

  {\rm (a)} For each $n \ge 1$, we have
  \begin{equation*}
    \Trace(\rho_{n+1}^V(b_n)) = \sum_{i_1, \cdots, i_n} \Trace_V^{}(s_{i_1} \cdots s_{i_n} t_{i_n} \cdots t_{i_1}).
  \end{equation*}

  {\rm (b)} If $A$ is involutive, $\Trace(\rho_{n+1}^V(b_n)) = \Trace_V^{}(u^{-n})$.
\end{lemma}

\begin{proof}
  (a) Let $f_0, \cdots, f_n \in \End(V)$ be linear endomorphisms on $V$. Define a linear map $f: V^{\otimes(n+1)} \to V^{\otimes(n+1)}$ by
  \begin{equation*}
    f(v_0 \otimes \cdots \otimes v_n) = f_1(v_1) \otimes \cdots \otimes f_n(v_n) \otimes f_0(v_0)
  \end{equation*}
  for all $v_0, \cdots, v_n \in V$. Then we have $\Trace(f) = \Trace(f_1 \circ \cdots \circ f_n \circ f_0)$ by direct calculation. Applying this formula to $\rho_{n+1}^V(b_n)$ the assertion follows.

  (b) Since $A$ is involutive, $u$ is central in $A$ and its inverse is given by $u^{-1} = \sum_i t_i s_i$ (see Proposition~\ref{prop:qt-hopf}). Thus we have
  \begin{equation*}
    \sum_{i_1, \cdots, i_n} t_{i_1} \cdots t_{i_n} s_{i_n} \cdots s_{i_1} = u^{-n}.
  \end{equation*}
  This implies $\Trace(\rho_{n+1}^V(b_n)) = \Trace_V^{}(u^{-n})$.
\end{proof}

\begin{remark}
  We can avoid the large part of the calculation. The proof will be much easier if we use Kauffman's beads arguments \cite{MR1734410} with suitable modification.
\end{remark}

The order of the antipode of a finite-dimensional Hopf algebra $H$ equals to the order of the antipode of $\mathcal{D}(H)$. The following description of $\tau(b_n; H)$ is a direct consequence of Lemma \ref{lem:trace-1}.

\begin{lemma}
  \label{lem:trace-2}
  For a finite-dimensional involutive Hopf algebra $H$, we have
  \begin{equation*}
    \tau(b_n; H) = \Trace_{\mathcal{D}(H)}^{}(u^{-n})
  \end{equation*}
  where $u$ is the Drinfeld element of the Drinfeld double $\mathcal{D}(H)$.
\end{lemma}

\begin{remark}
  \label{rem:FS-indicator}
  Similarly, we have $\tau(b_n^{-1}; H) = \Trace_{\mathcal{D}(H)}^{}(u^n)$ under the same assumption as the above lemma. If the characteristic of $k$ is zero, $\dim(H)^{-1} \tau(b_n^{-1}; H)$ equals to $\nu_n(\Trace_H^{})$ where $\nu_n$ is the $n$-th Frobenius-Schur indicator \cite{MR2213320}. The following lemma is only a well-known property of higher Frobenius-Schur indicators.
\end{remark}

\begin{lemma}
  \label{lem:trace-3}
  Let $G$ be a finite group. For each positive number $n$, we have
  \begin{equation*}
    \tau(b_n; kG) = |G| \cdot \# \{ g \in G \mid g^n = 1 \}.
  \end{equation*}
\end{lemma}

\begin{proof}
  Let $u$ be the Drinfeld element of $\mathcal{D}(kG)$. We calculate $\Trace_{\mathcal{D}(kG)}(u^{-n})$ in view of Lemma \ref{lem:trace-2}. We have $u^{-n} = \sum_{g \in G} e_g \bicross g^n$ by induction on $n$. By the definition of the multiplication, we have
  \begin{equation*}
    \Trace_{\mathcal{D}(kG)}^{}(e_x \bicross g) = \delta_{g, 1} |G|
  \end{equation*}
  for all $x, g \in G$. Thus, we have
  \begin{equation*}
    \tau(b_n; kG) = \Trace_{\mathcal{D}(kG)}(u^{-n}) = |G| \cdot \# \{ g \in G \mid g^n = 1 \}.
  \end{equation*}
\end{proof}

\begin{proof}[Proof of Theorem~\ref{thm:main-theorem}]
  Let $o_n(G)$ be the number of elements of order $n$ in $G$. By Lemma~\ref{lem:trace-3}, we have
  \begin{equation*}
    \frac{1}{|G|} \tau(b_n; \mathbb{C}G) = \sum_{d | n} o_d(G)
  \end{equation*}
  where the sum is taken over all positive integer $d$ that divides $n$. Applying the M{\"o}bius inversion formula to this equation, we have
  \begin{equation*}
    o_n(G) = \frac{1}{|G|} \sum_{d | n} \mu\left(\frac{n}{d}\right) \tau(b_d; \mathbb{C}G)
  \end{equation*}
  where $\mu$ is the M{\"o}bius function. If $G$ and $G'$ are $k$-isocategorical finite groups, then $\tau(b_d; \mathbb{C}G) = \tau(b_d; \mathbb{C}G')$ by Theorem~\ref{thm:base-change-theorem} (b). Hence we have $o_n(G) = o_n(G')$.
\end{proof}

Finally, we give some remarks on monoidal Morita invariants $\tau(b_n; -)$ and $\tau(b_n^{-1}; -)$. Until the end of this section, the base field $k$ is assumed to be an algebraically closed field of characteristic zero. For a finite-dimensional semisimple Hopf algebra $H$ and an integer $n$, we set
\begin{equation*}
  \omega_n(H) = \frac{1}{\dim(H)} \Trace_{\mathcal{D}(H)}^{}(u^n)
\end{equation*}
where $u$ is the Drinfeld element of the Drinfeld double $(\mathcal{D}(H), \mathcal{R})$. This is a monoidal Morita invariant by Lemma~\ref{lem:trace-2} and Remark~\ref{rem:FS-indicator}.

Let $V$ be a finite-dimensional $H$-module with character $\chi$. For a positive integer $n$, the {\em $n$-th Frobenius-Schur indicator} of $\chi$ is the number
\begin{equation*}
  \nu_n(\chi) := \sum \chi(\Lambda_{(1)} \Lambda_{(2)} \cdots \Lambda_{(n)})
\end{equation*}
where $\Lambda \in H$ is the integral such that $\varepsilon(\Lambda) = 1$. If $\widetilde{\chi}$ is the character of the induced module $\mathcal{D}(H) \otimes_H V$, $\nu_n(\chi) = \dim(H)^{-1} \, \widetilde{\chi}(u^n)$ \cite{MR2213320}. In particular, $\omega_n(H) = \nu_n(\chi_H^{})$ where $\chi_H^{}$ is the character of the regular representation $H$. The following lemma is essentially given in \cite{MR2213320}.

\begin{lemma}
  \label{lem:trace-4}
  $\omega_r(H) = 1$ if $r$ is coprime to $\dim(H)$. In particular, $\omega_{\pm 1}(H) = 1$.
\end{lemma}
\begin{proof}
  We first prove the case when $r = 1$. Let $\lambda \in H^*$ be the integral on $H$ such that $\langle \lambda, \Lambda \rangle = 1$. Then we have $\langle \lambda, 1 \rangle = \dim(H)$, and hence $\chi_H^{} = \lambda$ (see \cite[Chapter~7]{MR1786197}). Therefore, $\omega_1(H) = \nu_1(\chi_H^{}) = \langle \lambda, \Lambda \rangle = 1$.

  Next, we prove the general case. Since $H$ is semisimple, $u^{\dim(H)^3}$ = 1 \cite[Theorem~4]{MR1689203}. Let $\{ V_i \}_{i \in I}$ be representatives of isomorphism classes of simple $\mathcal{D}(H)$-modules. We denote $\dim(V_i)$ by $d_i$. Since $u$ is central, $u$ acts on $V_i$ as scalar. We denote this scalar by $u_i$. $u_i$'s are $\dim(H)^3$-th root of unity.
  By Artin-Wedderburn theorem, we have an isomorphism $\mathcal{D}(H) \cong \oplus_{i \in I} V_i^{\oplus d_i}$ of left $\mathcal{D}(H)$-modules. This decomposition yields
  \begin{equation*}
    \Trace_{\mathcal{D}(H)}(u^n) = \sum_{i \in I} u_i^n d_i^2
  \end{equation*}
  for any integer $n$.

  Let $\zeta$ be a primitive $\dim(H)^3$-th root of unity. Since $r$ is coprime to $\dim(H)$, there exists $\sigma \in {\rm Gal}(\mathbb{Q}[\zeta]/\mathbb{Q})$ such that $\sigma(\zeta) = \zeta^r$. Then, $\sigma(u_i) = u_i^r$ for all $i \in I$. Thus,
  \begin{equation*}
    \Trace_{\mathcal{D}(H)}^{}(u^r) = \sum_{i \in I} u_i^r d_i^2
    = \sum_{i \in I} \sigma(u_i) d_i^2 = \sigma(\dim(H)) = \dim(H).
  \end{equation*}
  This implies $\omega_r(H) = 1$.
\end{proof}

\section{Relations to low-dimensional topology}

\label{sec:topology}

We discuss relations between our invariants and the construction of invariants of closed 3-manifolds due to Reshetikhin and Turaev. Throughout this section, the base field $k$ is assumed to be algebraically closed of characteristic zero.

\subsection{Ribbon categories}

A {\em ribbon category} is a braided monoidal category with left duality and balancing isomorphism (see \cite[Chapter~I]{MR1292673} and \cite[Chapter~XIV]{MR1321145}). The balancing isomorphism is denoted by $\theta_V^{}: V \to V$. In this section, we use the {\em graphical calculus} (\cite[I.1.6]{MR1292673}, \cite[Chapter XIV]{MR1321145}) which is a pictorial technique to represent morphisms in ribbon categories.

Let $\mathcal{C}$ be a ribbon category. We say that an oriented framed link is $\mathcal{C}$-colored if each of its component is labeled with an object of $\mathcal{C}$. Every $\mathcal{C}$-colored oriented framed link $L$ defines a morphism $\mathbf{1} \to \mathbf{1}$ in $\mathcal{C}$ (see \cite[Chapter~I]{MR1292673}). This morphism is called the {\em operator valued invariant} of $L$.

Let $(A, R)$ be a quasitriangular Hopf algebra. A central invertible element $\theta \in A$ is called a {\em ribbon element} if we have $S(\theta) = \theta$ and $\Delta(\theta) = (R_{21}R)(\theta \otimes \theta)$. A {\em ribbon Hopf algebra} is a quasitriangular Hopf algebra equipped with a ribbon element. If $(A, R, \theta)$ is a ribbon Hopf algebra, finite-dimensional left $A$-modules form a ribbon category with braiding $c^R$ and a balancing isomorphism given by $\theta_V: V \to V$, $v \mapsto \theta \cdot v$ ($v \in V$). We denote this ribbon category by $\finMod(A, R, \theta)$, or simply by $\finMod(A)$ if $R$ and $\theta$ are obvious.

The quantum trace \cite[I.1.5]{MR1292673} of an endomorphism $f: V \to V$ in a ribbon category is denoted by $\trace_q(f)$. For an object $V$ of a ribbon category, the quantum dimension $\dim_q(V)$ is defined by $\dim_q(V) = \trace_q(\id_V)$. We argue some properties of the quantum trace in $\finMod(A)$.

\begin{lemma}
  \label{lem:quantum-trace}
  Notations are as above.

  {\rm (a)} Let $f: V \to V$ be a morphism in $\finMod(A)$. Then, we have
  \begin{equation*}
    \trace_q(f) = \Trace(v \mapsto u \theta \cdot f(v))
  \end{equation*}
  where $u$ is the Drinfeld element of $(A, R)$. Here, $\Trace$ means the usual trace of linear endomorphisms.

  {\rm (b)} Consider the commutative diagram of morphisms in $\finMod(A)$
  \begin{equation*}
    \begin{CD}
      0 @>>> V' @>{i}>> V @>{p}>> V'' @>>> 0 \\
      @. @V{f'}VV @V{f}VV @VV{f''}V \\
      0 @>>> V' @>>{i}> V @>>{p}> V'' @>>> 0
    \end{CD}
  \end{equation*}
  with exact rows. We have $\trace_q(f) = \trace_q(f') + \trace_q(f'')$.
\end{lemma}

\begin{proof}
  (a) is well known, see for example \cite[Proposition~XIV.6.4]{MR1321145}. (b) is obvious by (a) and usual properties of the trace.
\end{proof}

Let $\mathcal{C}$ be a ribbon category, $L$ an oriented framed link with components $L_1, \cdots, L_m$. If we color each component $L_i$ with an object $V_i \in \mathcal{C}$, we obtain an element $\End_{\mathcal{C}}(\mathbf{1})$ via the operator-valued invariant. We denote this element by $\langle L; V_1, \cdots, V_m \rangle$.

\begin{lemma}
  \label{lem:multiadditivity}
  Let $L$ be an oriented framed link with $m$-components. $\langle L; -, \cdots, - \rangle$ is multiadditive in the following sense: Let $V_1, \cdots, V_m \in \finMod(A)$. If there exists an exact sequence
  \begin{equation*}
    \begin{CD} 0 @>>> V_i' @>>> V_i @>>> V_i'' @>>> 0 \end{CD}
  \end{equation*}
  in $\finMod(A)$ for some $i$, we have
  \begin{equation*}
    \begin{split}
      \langle L; V_1, \cdots, V_m \rangle
      & = \langle L; V_1, \cdots, V_{i-1}, V_i', V_{i+1}, \cdots, V_m \rangle \\
      & + \langle L; V_1, \cdots, V_{i-1}, V_i'', V_{i+1}, \cdots, V_m \rangle.
    \end{split}
  \end{equation*}
\end{lemma}

\begin{proof}
  First, color each component $L_j$ with $V_j$ except for the $i$-th one. Then cut the $i$-th component of $L$ and form a partially colored ribbon graph $T$ so that we can obtain $L$ by closing $T$. For each $V \in \finMod(A)$, we obtain a morphism $\eta_V^{}: V \to V$ in $\finMod(A)$ by coloring the uncolored component with $V$. By graphical calculus, we have
  \begin{equation*}
    \trace_q(\eta_V^{}) = \langle L; V_1, \cdots, V_{i-1}, V, V_{i+1}, \cdots, V_m \rangle.
  \end{equation*}
  The family $\eta$ is a natural morphism in $\finMod(A)$. Thus, we have
  \begin{equation*}
    \trace_q(\eta_{V_i}^{}) = \trace_q(\eta_{V_i'}^{}) + \trace_q(\eta_{V_i''}^{})
  \end{equation*}
  by Lemma~\ref{lem:quantum-trace} (b). This completes the proof.
\end{proof}

The following lemma is due to Etingof and Gelaki \cite{MR1617921}.

\begin{lemma}
  \label{lem:ribbon-Hopf}
  Let $(A, R)$ be a quasitriangular Hopf algebra with $u$ the Drinfeld element. If $A$ is finite-dimensional semisimple and cosemisimple, $u^{-1}$ is a ribbon element of $(A, R)$.
\end{lemma}

In the ribbon category $\finMod(A, R, u^{-1})$, the quantum trace and the quantum dimension reduce to the usual trace and the dimension by Lemma~\ref{lem:quantum-trace}~(a).

\subsection{Reshetikhin-Turaev invariant of closed 3-manifolds}

By a ``manifold'' we mean an oriented connected topological manifold. A {\em modular category} \cite[Definition~3.1.1]{MR1797619} is a semisimple ribbon category with a finite number of simple objects satisfying a certain non-degeneracy condition. Reshetikhin and Turaev \cite{MR1036112} introduced a method of constructing an invariant of closed 3-manifold using a modular category. Let us briefly describe their construction following \cite{MR1292673} and \cite{MR1797619}.

Let $\mathcal{C}$ be a modular category with $\{ V_i \}_{i \in I}$ representatives of isomorphism classes of simple objects of $\mathcal{C}$. We denote by $d_i$ the quantum dimension of $V_i$. Since $\End_{\mathcal{C}}(V_i) = k$ by definition, we can define $\theta_i \in k$ by $\theta_{V_i} = \theta_i \cdot \id_{V_i}$ for each $i \in I$. Set $p^{\pm} = \sum_{i \in I} \theta_i^{\pm 1} d_i^2$ and $D = \sqrt{p^+ p^-}$. Then the numbers $p^\pm$ and $D$ are nonzero \cite[Theorem~3.1.7]{MR1797619}.

For a framed link $L$ with components $L_1, \cdots, L_m$, we fix an arbitrary orientation of $L$ and set
\begin{equation*}
  \label{eq:framed-link-eval}
  \{ L \} = \sum_{i_1, \cdots, i_m \in I} \langle L; V_{i_1}, \cdots, V_{i_m} \rangle \, d_{i_1} \cdots d_{i_m}.
\end{equation*}
The right hand side does not depend on the numbering of components and the choice of orientation of $L$.

Now we describe the {\em Reshetikhin-Turaev invariant} ${\bf RT}_{\mathcal{C}}^{}$ of a closed 3-manifold associated with $\mathcal{C}$. Let $M$ be a closed 3-manifold. By a classical result, any closed 3-manifold can be obtained by the so-called {\em Dehn surgery} on $S^3$ along a certain framed link. Fix a framed link $L$ yielding $M$. Then, ${\bf RT}_{\mathcal{C}}^{}(M)$ is given by
\begin{equation*}
  \label{eq:RT-invariant}
  {\bf RT}_{\mathcal{C}}^{}(M) = D^{-|L|-1} \left(\frac{p^+}{p^-}\right)^{\frac{1}{2}\sigma(L)} \{ L \}
\end{equation*}
where $|L|$ is the number of components of $L$ and $\sigma(L)$ is the so-called {\em wreath number} of $L$ (see \cite[II.2.1]{MR1292673} for its definition). The right hand side does not depend on the choice of $L$ and thus ${\bf RT}_{\mathcal{C}}^{}(M)$ is an invariant of the closed 3-manifold $M$.

\begin{remark}
  We defined $D$ to be $\sqrt{p^+p^-}$. This exists since we work over an algebraically closed field of characteristic zero. The definition of ${\bf RT}_{\mathcal{C}}^{}$ depends on the choice of $D$, that is, the choice of square roots of $p^+p^-$. Therefore, we need to fix $D$ to define ${\bf RT}_{\mathcal{C}}^{}$.
\end{remark}

\subsection{Modular categories arising from Hopf algebras}

Let $H$ be a finite-dimensional semisimple Hopf algebra and $u$ the Drinfeld element for the Drinfeld double $\mathcal{D}(H)$. Then $\finMod(\mathcal{D}(H), \mathcal{R}, u^{-1})$ is a modular category \cite[Lemma~1.1]{MR1617921}. We denote by ${\bf RT}_{\mathcal{D}(H)}^{}$ the Reshetikhin-Turaev invariant associated to this modular category.

Let us describe the invariant ${\bf RT}_{\mathcal{D}(H)}$. First, we compute numbers $p^{\pm}$ and $D$. Let $\{ V_i \}_{i \in I}$ be representatives of isomorphism classes of irreducible $\mathcal{D}(H)$-modules. By Artin-Wedderburn theorem, we have an isomorphism $\mathcal{D}(H) \cong \bigoplus_{i \in I} V_i^{\oplus d_i}$ of left $\mathcal{D}(H)$-modules. Note that $\theta_i$ is the unique eigenvalue of the action of central element $u^{-1}$ on $V_i$. By Lemma~\ref{lem:trace-4}, we have
\begin{equation*}
  p^{\pm} = \sum_{i \in I} \theta_i^{\pm 1} d_i^2 = \Trace_{\mathcal{D}(H)}(u^{\mp 1}) = \dim(H).
\end{equation*}
This allows us to choose $D = \sqrt{p^+p^-}$ to be $\dim(H)$.

\begin{remark}
  $\zeta = (p^+/p^-)^{1/6}$ is known to be a root of unity in general. When the base field is $\mathbb{C}$, we can write $\zeta = \exp(2 \pi c \sqrt{-1} / 24)$ for some $c$. $c$ is called the {\em central charge} for the theory \cite[Remark~3.1.20]{MR1797619}. In our cases, we have $p^\pm = \dim(H)$ as described above. This implies that the central charge of $\finMod(\mathcal{D}(H))$ is zero.
\end{remark}

Next, we argue $\{ L \}$ where $L$ is a framed link with $m$ components. Note that every $d_i$ is a positive integer. By Lemma~\ref{lem:multiadditivity},
\begin{equation*}
  \{ L \} = \sum_{i_1, \cdots, i_m} \left\langle L; V_{i_1}^{\oplus d_{i_1}}, \cdots, V_{i_m}^{\oplus d_{i_m}} \right\rangle = \langle L; \mathcal{D}(H), \cdots, \mathcal{D}(H) \rangle.
\end{equation*}
Summarizing, we have the following theorem.

\begin{theorem}
  \label{thm:topology-1}
  The Reshetikhin-Turaev invariant ${\bf RT}_{\mathcal{D}(H)}^{}$ is given as follows: If a closed 3-manifold $M$ is obtained by surgery on $S^3$ along a framed link $L$,
  \begin{equation*}
    {\bf RT}_{\mathcal{D}(H)}^{}(M) = \dim(H)^{-|L|-1} \langle L; \mathcal{D}(H), \cdots, \mathcal{D}(H) \rangle.
  \end{equation*}
\end{theorem}
 
For a braid $b \in B_n$, we denote by $\widehat{b}$ the framed link obtained by closing $b$. A graphical calculus gives the equation
\begin{equation*}
  \Trace(\rho_n^{\mathcal{D}(H)}(b)) = \langle \,\widehat{b}\,; \mathcal{D}(H), \cdots, \mathcal{D}(H) \rangle.
\end{equation*}
Thus, if a closed 3-manifold $M$ is obtained by surgery along $\widehat{b}$, we have
\begin{equation*}
  {\bf RT}_{\mathcal{D}(H)}^{}(M) = \dim(H)^{- |\,\widehat{b}\,| - 1} \, \tau(b; H).
\end{equation*}
Note that any framed link can be obtained by closing a certain braid (for ordinary links, this fact is known as Alexander's theorem). Thus, any closed 3-manifold can be obtained by surgery along $\widehat{b}$ for some braid $b$. Let $H$ and $L$ be finite-dimensional semisimple Hopf algebras. If $H$ and $L$ are monoidally Morita equivalent, we have $\dim(H) = \dim(L)$ by Lemma~\ref{lem:key-lemma}. Summarizing, we have the following theorem.

\begin{theorem}
  \label{thm:topology-2}
  $H$ and $L$ are as above. Then we have ${\bf RT}_{\mathcal{D}(H)}^{} = {\bf RT}_{\mathcal{D}(L)}^{}$.
\end{theorem}

Let $G$ be a finite group and $\omega: G \times G \times G \to \mathbb{C}^\times$ a normalized 3-cocycle. Dijkgraaf and Witten \cite{MR1048699} introduced a method of constructing an invariant of closed 3-manifolds using a pair $(G, \omega)$ (see also \cite{MR1192735}). Let us denote this invariant by $Z_{G, \omega}$. When $\omega$ is the trivial 3-cocycle, by definition, we have
\begin{equation*}
  Z_{G, 1}(M) = \frac{1}{|G|} \# \Hom(\pi_1(M), G)
\end{equation*}
where $\pi_1(M)$ is the fundamental group of $M$.

On the other hand, Alts\"uler and Coste \cite{MR1188498} introduced a method of constructing an invariant of 3-manifolds using the modular category of finite-dimensional modules over the quasi-Hopf algebra $\mathcal{D}^\omega(G)$, which is defined to be a certain deformation of $\mathcal{D}(\mathbb{C}G)$. When $\omega$ is the trivial 3-cocycle, this invariant is equal to ${\bf RT}_{\mathcal{D}(\mathbb{C}G)}$. Alts\"uler and Coste conjectured in \cite{MR1188498} and Sato and Wakui proved in \cite[Corollary 5.5]{MR2002616} that Alts\"uler-Coste invariant is equal to the Dijkgraaf-Witten invariant $Z_{G, \omega}$. Thus we have
\begin{equation*}
  {\bf RT}_{\mathcal{D}(\mathbb{C}G)}^{}(M) = Z_{G, 1}(M) = \frac{1}{|G|} \# \Hom(\pi_1(M), G)
\end{equation*}
for all closed 3-manifold $M$. This gives rise to the following theorem.

\begin{theorem}
  \label{thm:topology-3}
  Let $b$ be a braid, $m$ the number of components of $\widehat{b}$ and $M$ the closed 3-manifold obtained by surgery along $\widehat{b}$. Then we have
  \begin{equation*}
    \tau(b; \mathbb{C}G) = |G|^{m} \, \# \Hom(\pi_1(M), G)
  \end{equation*}
\end{theorem}

Combining Theorem~\ref{thm:base-change-theorem} and this theorem, we have Theorem~\ref{thm:main-theorem-2}. Note that Lemma~\ref{lem:trace-3} is a special case of the above theorem when $M$ is the lens space $L(n, 1)$ whose fundamental group is the cyclic group of order $n$. In fact, it is well known in the theory of surgery that framed link $\widehat{b}_n$ yields $L(n, 1)$. ($\widehat{b}_n$ is isotopic to the trivial knot with the framing $+n$.)

\section{Examples}

\label{sec:examples}

\subsection{Numbers of homomorphisms from quaternion groups}

For an integer $m \ge 2$, set $Q_{4m} = \langle x, y \mid x^{2m} = 1, y^2 = x^m, yxy^{-1} = x^{-1} \rangle$. $Q_{8}$ is the quaternion group. In general, $Q_{4m}$ is called the generalized quaternion group of order $4m$. The following theorem is an example of our monoidal Morita invariants.

\begin{theorem}
  \label{thm:hom-Q8}
  Set $b = \sigma_1^4 \in B_2$. Then, for any finite group $G$, we have
  \begin{equation*}
    \tau(b; \mathbb{C}G) = |G|^2 \, \# \Hom(Q_8, G).
  \end{equation*}
\end{theorem}
\begin{proof}
  For simplicity, we write $\mathcal{D}(\mathbb{C}G)$ by $\mathcal{D}(G)$. Let $\mathcal{R}$ be the universal $R$-matrix of $\mathcal{D}(G)$. Then, in a similar way to the proof of Lemma~\ref{lem:trace-3}, we have
  \begin{equation*}
    \begin{split}
      \tau(b; \mathbb{C}G)
      &= \Trace_{\mathcal{D}(G) \otimes \mathcal{D}(G)}^{}
      \left( \left( \mathcal{R}_{21} \mathcal{R} \right)^2 \right) \\
      &= \sum_{g, h \in G} \Trace_{\mathcal{D}(G)}\left( e_g \bicross g^{-1} h g h \right)
      \Trace_{\mathcal{D}(G)}\left( e_{g h g^{-1}} \bicross g h^{-1} g h \right) \\
      &= |G|^2 \, \# \{ (g, h) \in G \times G \mid \mbox{$g^{-1} h g h = 1$, $g h^{-1} g h = 1$} \} \\
      &= |G|^2 \, \# \Hom(Q, G)
    \end{split}
  \end{equation*}
  where $Q$ is the group defined by generators $g$ and $h$ with relations $g^{-1} h g h = g h^{-1} g h = 1$. $Q$ is isomorphic to $Q_8$ via a map $Q_8 \to Q$ given by $x \mapsto g$, $y \mapsto h$.
\end{proof}

\begin{remark}
  The list of all finite subgroups of $SO(4)$ which can act freely on $S^3$ is known (see, e.g., \cite[\S 6]{MR0426001}) and contains $Q_{4m}$ for all $m \ge 2$. If $\Gamma$ is such a finite group, the quotient space $S^3/\Gamma$ is an orientable closed 3-manifold with fundamental group $\Gamma$ (spherical manifolds). Thus, in view of Theorem~\ref{thm:topology-3}, there exists a braid $b$ such that $\tau(b; \mathbb{C}G) = |G|^m \,\#\Hom(\Gamma, G)$ where $m$ is the number of components of $\widehat{b}$. The above theorem may be considered as a special case of this fact. However, the author does not know whether the closed 3-manifold $S^3/Q_8$ is obtained by surgery along $\widehat{b}$ with $b = \sigma_1^4$.
\end{remark}

\subsection{Categorical rigidity of finite groups of small order}

In this subsection, we argue categorical rigidity of finite groups of small order and prove the following theorem.

\begin{theorem}
  All finite groups of orders less than 32 are categorically rigid.
\end{theorem}

\begin{proof}
  We first note that all finite abelian groups are categorically rigid over an arbitrary field. In fact, if $H$ is a finite-dimensional commutative Hopf algebra, every finite-dimensional Hopf algebra which is monoidally Morita equivalent to $H$ is isomorphic to $H$ (see \cite[Remark~3.8]{MR1408508}). Therefore, we consider non-abelian finite groups.

  By using Theorem~\ref{thm:main-theorem} and the above-mentioned fact, we can conclude that all finite groups of orders less than $32$ except for $16$ are categorically rigid over an arbitrary field. However, this theorem is not sufficient to give the complete classification of groups of order $16$.

   Let us argue categorical rigidity of groups of order 16. As is well known, there are exactly nine non-abelian groups of order $16$ up to isomorphism. Using Theorem~\ref{thm:main-theorem}, we conclude that five of them are categorically rigid over an arbitrary field. The rest of them consists of two pairs for which Theorem~\ref{thm:main-theorem} fails to work. The first pair consists of
  \begin{equation*}
    G_1 = Q_8 \times \mathbb{Z}_2 \quad \text{and} \quad
    G_2 = \langle x, y \mid x^4 = y^4 = 1, y x y^{-1} = x^{-1} \rangle.
  \end{equation*}
  They are not isomorphic, but, for each positive integer $n$, the number of elements of order $n$ in them are equal. We conclude that $G_1$ and $G_2$ are not isocategorical by using Theorem~\ref{thm:hom-Q8}. In fact, we have
  \begin{equation*}
    \# \Hom(Q_8, G_1) = 112 \quad \text{and} \quad
    \# \Hom(Q_8, G_2) = 16.
  \end{equation*}

  The second pair $(F_1, F_2)$ is given as follows. Set
  \begin{equation*}
    F = \langle x, y \mid x^4 = y^2 = 1, xy = yx \rangle \quad \text{and} \quad
    C_2 = \langle s \mid s^2 = 1 \rangle,
  \end{equation*}
  and define automorphisms $f_1$ and $f_2$ on $F$ respectively by
  \begin{equation*}
    f_1(x) = x,  \quad f_1(y) = x^2 y \quad \text{and} \quad
    f_2(x) = xy, \quad f_2(y) = y.
  \end{equation*}
  $F_1$ and $F_2$ are semidirect products $F \rtimes C_2$ where $s \in C_2$ acts on $F$ respectively by $f_1$ and $f_2$. They are not isomorphic since their abelianizations are different:
  \begin{equation*}
    F_1^{\rm ab} \cong \mathbb{Z}_2 \oplus \mathbb{Z}_2 \oplus \mathbb{Z}_2 \quad \text{and} \quad
    F_2^{\rm ab} \cong \mathbb{Z}_2 \oplus \mathbb{Z}_4.
  \end{equation*}
  This pair is an example for which both Theorem~\ref{thm:main-theorem} and Theorem~\ref{thm:hom-Q8} fail to work. In fact, they have an equal number of elements of order $n$ for each positive integer $n$, and, moreover,
  \begin{equation*}
    \# \Hom(Q_8, F_1) = \# \Hom(Q_8, F_2) = 64.
  \end{equation*}

  If $k$ is an algebraically closed field of ${\rm char}(k) \ne 2$, $F_1$ and $F_2$ are not $k$-isocategorical since their Grothendieck rings are different. If ${\rm char}(k) = 2$, $kF_1$ and $kF_2$ are not even Morita equivalent. This can be proved as follows. First, we note that every irreducible representation of a finite $p$-group in characteristic $p$ is isomorphic to the trivial one. Therefore, if they are Morita equivalent, there is an isomorphism
  \begin{equation*}
    \Ext_{kF_1}^1(k, k) \cong \Ext_{kF_2}^1(k, k).
  \end{equation*}
  Now we recall that there is an isomorphism $\Ext_{kG}^1(k, k) \cong \Hom(G^{\rm ab}, k)$ for every group $G$. We have
  \begin{equation*}
    \Ext_{kF_1}^1(k, k) \cong k^3 \quad \text{and} \quad
    \Ext_{kF_2}^1(k, k) \cong k^2.
  \end{equation*}
  This is a contradiction. Therefore, all finite groups of order $16$ are categorically rigid over an arbitrary field.
\end{proof}

\section*{Acknowledgment}

The author is grateful to Professor Mitsuhiro Takeuchi for his valuable comments on a draft of this paper. This work is supported by Grant-in-Aid for JSPS Fellows.

\appendix
\def\thesection{\Alph{section}}

\section{Similarity between permutation matrices}

\label{sec:similarity}

We denote by $\mathfrak{S}_n$ the symmetric group of degree $n$. For $\sigma \in \mathfrak{S}_n$, we denote by $P_\sigma$ the $n \times n$ matrix whose $(i, j)$-entry is $\delta_{\sigma(i), j}$ where $\delta$ is Kronecker's delta. We prove the following theorem.

\begin{theorem}
  \label{thm:similarity-theorem}
  $P_\sigma$ and $P_\tau$ are similar if and only if $\sigma$ and $\tau$ are conjugate.
\end{theorem}

The ``if'' part is clear since the map $\sigma \mapsto P_\sigma$ is a group homomorphism. Let us prove the ``only if'' part. As a first step, we characterize linear automorphisms of a finite-dimensional vector space which are represented by permutation matrices in some basis. Let $A = k[X, X^{-1}]$ be the Laurent polynomial ring with an indeterminate $X$. For a vector space $V$ and a linear automorphism $P$ on $V$, we denote by ${}_P V$ the $A$-module with the underlying space $V$ and the action given by $X \cdot v = P(v)$ for $v \in V$. Set
\begin{equation*}
  M(n) = k[X, X^{-1}] /(X^n - 1) \quad (n = 1, 2, \cdots).
\end{equation*}

\begin{lemma}
  \label{lem:characterize-permutations}
  Let $V$ be a finite-dimensional vector space and $P$ be a linear automorphism on $V$. $P$ is represented by a permutation matrix in some basis if and only if the $A$-module ${}_P V$ is isomorphic to a direct sum of $M(i)$'s.
\end{lemma}
\begin{proof}
  The ``if'' part is clear since the action of $X$ on $M(i)$ is represented by a permutation matrix in basis $\{ 1, X, \cdots, X^{i-1} \}$. We prove the ``only if'' part. Set $n = \dim(V)$. Assume that $P$ is represented by the permutation matrix $P_\sigma$ for some $\sigma \in \mathfrak{S}_n$ in basis $\{ e_1, \cdots, e_n \}$. Let
  \begin{equation}
    \label{eq:sigma-orbit-decomp}
    \{ 1, 2, \cdots, n \} = \mathcal{O}_1 \sqcup \cdots \sqcup \mathcal{O}_r
  \end{equation}
  be the $\sigma$-orbit decomposition. $V_i = {\rm span}_k \{ e_s \mid s \in \mathcal{O}_i \}$ is an $A$-submodule of ${}_P V$ isomorphic to $M(\# \mathcal{O}_i)$. Thus, we have an isomorphism
  \begin{equation*}
    V = V_1 \oplus \cdots \oplus V_r \cong M(\# \mathcal{O}_1) \oplus \cdots \oplus M(\# \mathcal{O}_r)
  \end{equation*}
  of $A$-modules.
\end{proof}

Actually, the $\sigma$-orbit decomposition (\ref{eq:sigma-orbit-decomp}) gives a cycle decomposition of the permutation $\sigma$. We say that a finite-dimensional $A$-module $M$ {\em admits a cycle decomposition} if $M$ is isomorphic to a direct sum of $M(i)$'s.

Set $V = k^n$. Note that two invertible $n \times n$-matrices $P$ and $Q$ are similar if and only if ${}_P V$ and ${}_Q V$ are isomorphic as $A$-modules. As we observed above, if $P = P_\sigma$ ($\sigma \in \mathfrak{S}_n$) is a permutation matrix, ${}_P V$ admits a cycle decomposition
\begin{equation*}
  {}_P V \cong M(1)^{\oplus c_1(\sigma)} \oplus \cdots \oplus M(n)^{\oplus c_n(\sigma)}
\end{equation*}
where $c_r(\sigma)$ is the number of cyclic permutations of length $r$ which appear in the cyclic decomposition of $\sigma$. If another permutation matrix $Q = P_\tau$ ($\tau \in \mathfrak{S}_n$) is similar to $P$, ${}_P V$ admits another cycle decomposition
\begin{equation*}
  {}_P V \cong M(1)^{\oplus c_1(\tau)} \oplus \cdots \oplus M(n)^{\oplus c_n(\tau)}.
\end{equation*}
Recall that $\sigma$ and $\tau$ are conjugate if and only if their cycle shape are same, i.e., $c_r(\sigma) = c_r(\tau)$ for all $r$. Theorem~\ref{thm:similarity-theorem} turns into the following statement: If
\begin{equation*}
  M(1)^{\oplus d_1} \oplus \cdots \oplus M(n)^{\oplus d_n} \cong
  M(1)^{\oplus e_1} \oplus \cdots \oplus M(n)^{\oplus e_n}
\end{equation*}
as $A$-modules, then $d_i = e_i$ for all $i$.

Note that $A$ has a Hopf algebra structure as a group algebra of an infinite cyclic group generated by $X \in A$. Thus, we can consider tensor products and dual modules of $A$-modules.

\begin{lemma}
  \label{lem:tensor-and-dual}
  {\rm (a)} $M(n) \otimes M(m) \cong M(mn/d)^{\oplus d}$ as an $A$-module where $d = \gcd(n, m)$ is the greatest common divisor of $n$ and $m$.

  {\rm (b)} $M(n)^* \cong M(n)$ as an $A$-module.
\end{lemma}
\begin{proof}
  {\rm (a)} The action of $X$ permutes $\{ X^i \otimes X^j \} \subset M(n) \otimes M(m)$. Easy combinatorial arguments completes the proof.

  {\rm (b)} Define $X_i^* \in M(n)^*$ by $X_i^*(X^j) = \delta_{i, n-j}$. Then the linear map $M(n) \to M(n)^*$ given by $X^i \mapsto X_i^*$ gives an isomorphism of $A$-modules.
\end{proof}

\begin{lemma}
  \label{lem:hom-dimension}
  $\dim \Hom_A(M(n), M(m)) = \gcd(n, m)$.
\end{lemma}
\begin{proof}
  First, we prove the case when $n = 1$. Let $f: M(1) \to M(m)$ be an $A$-linear map. If $f(1) = \sum_{i = 0}^{m-1} c_i X^i$ ($c_i \in k$), we have $c_0 = \cdots = c_{m-1}$ by $A$-linearity of $f$. Thus, we have $\dim \Hom_A(M(1), M(m)) = 1$.

  Now we prove the general case. By Lemma~\ref{lem:tensor-and-dual}, we have isomorphisms
  \begin{equation*}
    \begin{split}
      \Hom_A(M(n), M(m)) & \cong \Hom_A(M(1), M(m) \otimes M(n)^*) \\
      & \cong \Hom_A(M(1), M(l))^{\oplus d} \\
    \end{split}
  \end{equation*}
  where $d = \gcd(n, m)$ and $l = mn/d$. Therefore, we have
  \begin{equation*}
    \dim \Hom_A(M(n), M(m)) = d \cdot \dim \Hom_A(M(1), M(l)) = d.
  \end{equation*}
\end{proof}

\begin{proof}[Proof of Theorem~\ref{thm:similarity-theorem}]
  Suppose that
  \begin{equation*}
    V := M(1)^{\oplus e_1} \oplus \cdots \oplus M(n)^{\oplus e_n} \cong
    M(1)^{\oplus f_1} \oplus \cdots \oplus M(n)^{\oplus f_n}
  \end{equation*}
  as $A$-modules. Our aim is to prove $e_i = f_i$ for each $i$. Set $d_i = \dim \Hom_A(M(i), V)$. By Lemma~\ref{lem:hom-dimension}, we have equations
  \begin{equation*}
    d_i = \sum_{j = 1}^m \gcd(i, j) e_j
    \quad \mbox{and} \quad d_i = \sum_{j = 1}^m \gcd(i, j) f_j.
  \end{equation*}
  Thus, we have a linear equation
  \begin{equation*}
    \Phi_m \left( \begin{array}{c} e_1 \\ \vdots \\ e_n  \end{array} \right)
    = \Phi_m \left( \begin{array}{c} f_1 \\ \vdots \\ f_n \end{array} \right)
  \end{equation*}
  where $\Phi_m$ is an $m \times m$ matrix whose $(i, j)$-entry is $\gcd(i, j)$. The determinant of $\Phi_m$, which is known as Smith's determinant (see, e.g., \cite{MR1444107}), equals $\varphi(1) \cdots \varphi(n)$ where $\varphi$ is Euler's totient function. In particular, $\Phi_m$ is invertible, and thus we have $e_i = f_i$ for each $i$.
\end{proof}


\begin{thebibliography}{10}
\expandafter\ifx\csname url\endcsname\relax
  \def\url#1{\texttt{#1}}\fi
\expandafter\ifx\csname urlprefix\endcsname\relax\def\urlprefix{URL }\fi
\expandafter\ifx\csname href\endcsname\relax
  \def\href#1#2{#2} \def\path#1{#1}\fi

\bibitem{MR1810480}
P.~Etingof, S.~Gelaki, Isocategorical groups, Internat. Math. Res. Notices~(2)
  (2001) 59--76.

\bibitem{MR1292673}
V.~G. Turaev, Quantum invariants of knots and 3-manifolds, Vol.~18 of de
  Gruyter Studies in Mathematics, Walter de Gruyter \& Co., Berlin, 1994.

\bibitem{MR1408508}
P.~Schauenburg, Hopf bi-{G}alois extensions, Comm. Algebra 24~(12) (1996)
  3797--3825.

\bibitem{MR1797619}
B.~Bakalov, A.~Kirillov, Jr., Lectures on tensor categories and modular
  functors, Vol.~21 of University Lecture Series, American Mathematical
  Society, Providence, RI, 2001.

\bibitem{MR1786197}
S.~D{\u{a}}sc{\u{a}}lescu, C.~N{\u{a}}st{\u{a}}sescu, {\c{S}}.~Raianu, Hopf
  algebras, Vol. 235 of Monographs and Textbooks in Pure and Applied
  Mathematics, Marcel Dekker Inc., New York, 2001, an introduction.

\bibitem{MR1321145}
C.~Kassel, Quantum groups, Vol. 155 of Graduate Texts in Mathematics,
  Springer-Verlag, New York, 1995.

\bibitem{MR834465}
Y.~Doi, M.~Takeuchi, Cleft comodule algebras for a bialgebra, Comm. Algebra
  14~(5) (1986) 801--817.

\bibitem{MR1025154}
V.~G. Drinfeld, Almost cocommutative {H}opf algebras, Algebra i Analiz
  1~(2) (1989) 30--46.

\bibitem{MR1689203}
P.~Etingof, S.~Gelaki, On the exponent of finite-dimensional {H}opf algebras,
  Math. Res. Lett. 6~(2) (1999) 131--140.

\bibitem{MR1220770}
D.~E. Radford, Minimal quasitriangular {H}opf algebras, J. Algebra 157~(2)
  (1993) 285--315.

\bibitem{MR1734410}
L.~H. Kauffman, Right integrals and invariants of three-manifolds, in:
  Proceedings of the {K}irbyfest ({B}erkeley, {CA}, 1998), Vol.~2 of Geom.
  Topol. Monogr., Geom. Topol. Publ., Coventry, 1999, pp. 215--232
  (electronic).

\bibitem{MR2213320}
Y.~Kashina, Y.~Sommerh{\"a}user, Y.~Zhu, On higher {F}robenius-{S}chur
  indicators, Mem. Amer. Math. Soc. 181~(855) (2006) viii+65.

\bibitem{MR1617921}
P.~Etingof, S.~Gelaki, Some properties of finite-dimensional semisimple {H}opf
  algebras, Math. Res. Lett. 5~(1-2) (1998) 191--197.

\bibitem{MR1036112}
N.~Y. Reshetikhin, V.~G. Turaev, Ribbon graphs and their invariants derived
  from quantum groups, Comm. Math. Phys. 127~(1) (1990) 1--26.

\bibitem{MR1048699}
R.~Dijkgraaf, E.~Witten, Topological gauge theories and group cohomology, Comm.
  Math. Phys. 129~(2) (1990) 393--429.

\bibitem{MR1192735}
M.~Wakui, On {D}ijkgraaf-{W}itten invariant for {$3$}-manifolds, Osaka J. Math.
  29~(4) (1992) 675--696.

\bibitem{MR1188498}
D.~Altsch{\"u}ler, A.~Coste, Quasi-quantum groups, knots, three-manifolds, and
  topological field theory, Comm. Math. Phys. 150~(1) (1992) 83--107.

\bibitem{MR2002616}
  N.~Sato, M.~Wakui, {$(2+1)$}-dimensional topological quantum field theory with
  a {V}erlinde basis and {T}uraev-{V}iro-{O}cneanu invariants of 3-manifolds,
  in: Invariants of knots and 3-manifolds (Kyoto, 2001), Vol.~4 of Geom. Topol.
  Monogr., Geom. Topol. Publ., Coventry, 2002, pp. 281--294 (electronic).

\bibitem{MR0426001}
  P.~Orlik, Seifert manifolds, Lecture Notes in Mathematics, Vol. 291,
  Springer-Verlag, Berlin, 1972.

\bibitem{MR1444107}
  P.~Haukkanen, J.~Wang, J.~Sillanp{\"a}{\"a}, On {S}mith's determinant, Linear
  Algebra Appl. 258 (1997) 251--269.

\bibitem{MR2381536}
  Siu-Hung Ng and Peter Schauenburg.
  Higher {F}robenius-{S}chur indicators for pivotal categories.
  In {\em Hopf algebras and generalizations}, volume 441 of {\em Contemp. Math.}, pages 63--90. Amer. Math. Soc., Providence, RI, 2007.

\bibitem{MR2366965}
  Siu-Hung Ng and Peter Schauenburg.
  Central invariants and higher indicators for semisimple quasi-{H}opf algebras.
  {\em Trans. Amer. Math. Soc.}, 360(4):1839--1860, 2008.
\end{thebibliography}

\section*{Notes added in proof}

After this paper was accepted for publication, I received from Ng the following comments which show that Theorem~\ref{thm:main-theorem} follows from his results joint with Schauenburg under the assumption that the base field $k$ is an algebraically closed field of characteristic zero. I thank Ng for his valuable comments.

In \cite{MR2381536}, Ng and Schauenburg defined higher Frobenius-Schur indicators for pivotal monoidal categories and proved that these indicators are invariants under equivalences of such categories. Let $H$ and $L$ be finite-dimensional semisimple quasi-Hopf algebras. If $F: \Mod(H) \to \Mod(L)$ is an equivalence of monoidal categories, by \cite[Proposition~3.2]{MR2366965}, we have $\nu_m(V) = \nu_m(F(V))$ for every finite-dimensional $H$-module $V$. In particular, $\nu_m(H) = \nu_m(F(H)) = \nu_m(L)$. On the other hand, if $H = kG$ for some finite group $G$, we have
\begin{equation*}
  \nu_m(kG) = \# \{ g \in G \mid g^m = 1 \}.
\end{equation*}
This completes the proof.

\end{document}